\magnification=\magstephalf
\input amstex
\documentstyle{amsppt}
\pageheight{9truein}
\pagewidth{6.5truein}
 
\TagsOnRight
\loadbold
 
\define \gln {\text{\rm GL}_n}
\define \gl {\text{\rm GL}}
\define \gltwo {\text{\rm GL}_2}
\define \hofF {{\Cal H}(F)}
\define \hofG {{\Cal H}(G)}
\define \FL {F_{\Lambda}}
\define \fm {\frak m}
\define \fM {\Cal M}
\define \mofF {\fm (F)}
\define \MofF {\fM (F)}
\define \MofFL {\fM (\FL )}
\define \MofFLm {\fM (\FL ,m)}
\define \mofFL {\fm (\FL )}
\define \cont {\text{\rm cont}}
\define \disc {\text{\rm disc}}

\define \br {\Bbb R}
\define \be {\Bbb E}
\define \bc {\Bbb C}
\define \bz {\Bbb Z}
\define \bzp {\bz _p}
\define \bq {\Bbb Q}
\define \bqp {\bq _p}
\define \bqbar {\overline{\bq}}
\define \bqpbar {\overline{\bqp}}
\define \bcp {\bc _p}
\define \ux {\bold x}
\define \uX {\bold X}
\define \uY {\bold Y}
\define \uo {\bold 0}
\define \um {\bold m}
\define \ue {\bold e}
\define \ua {\bold a}
\define \uj {\bold j}
\define \uy {\bold y}
\define \uz {\bold z}
\define \uL {\bold L}
\define \uM {\bold M}
\define \uF {\bold F}
\define \uK {\bold K}
 
\topmatter
\title
Thue Equations and Lattices
\endtitle
 
\author
Jeffrey Lin Thunder
\endauthor
\address Department of Mathematics, Northern Illinois University,
DeKalb, IL 60115\endaddress
\email jthunder\@math.niu.edu\endemail
 
\abstract
We consider Diophantine equations of the kind $|F(x,y)|= m$,
where $F(X,Y )\in \bz [X,Y]$ is a homogeneous polynomial of degree $d\ge 3$ that has
non-zero discriminant and $m$ is a positive integer. We prove results that simplify
those of Stewart and provide heuristics for a conjecture of Stewart. 
\endabstract

\endtopmatter
 
\document
\baselineskip=20pt
 
\head Introduction\endhead
Suppose $F(X,Y)\in\bz [X,Y]$ is a homogeneous polynomial of degree $d\ge 3$ that has
non-zero discriminant and $m$ is a positive integer. In this paper we are concerned 
with the number of primitive, i.e., $x$ and $y$ are relatively prime, solutions
$(x,y)\in\bz ^2$ to the Thue equation
$$|F(x,y)|=m.\tag 1$$
Thue in [8] famously showed that the number of such solutions is necessarily finite under
the hypothesis that $F$ is irreducible over $\bq$. In fact,
his method enabled one to derive an upper bound on the number of such solutions; such an
upper bound would depend on $m$ and the polynomial $F$. Indeed, Lewis and Mahler in [4]
provided just such a bound. Their bound was an explicit function of $m$, $d$ and the height
of $F$. Previous to the result of Lewis and Mahler, Siegel had made the conjecture
that an upper bound could be obtained that was independent of the particular coefficients of
the polynomial $F$. Evertse proved this conjecture in his doctoral thesis (see [3]). 
A few years later
Bombieri and Schmidt [2] improved markedly on Evertse's bound, showing that the number of
primitive solutions to (1) is no more than some fixed (absolute) constant multiple of
$d^{1+\omega (m)}$, where $\omega (m)$ denotes the number of distinct prime factors of
$m$, as usual. Later Schmidt posited (see [6, chap. 3, Conjecture])
that the number of primitive solutions to (1) should
be bounded above by some multiple (possibly depending on $F$) of a power of $\log m$ when
$m>1$. 

Two years after the publication of Bombieri and Schmidt's result, Stewart [7] provided a
bound that was often (depending on the prime factorization of the parameter $m$) much
stronger than the bound of Bombieri and Schmidt. Stewart's main result was somewhat involved
and complicated to state, but one can easily state the following consequence. In what
follows, $D(F)$ denotes the discriminant of the form $F$.

\proclaim{Theorem (Stewart)} Suppose $F(X,Y)\in\bz [X,Y]$ is a homogeneous polynomial
of degree $d\ge 3$ with non-zero discriminant and content 1. Let $\epsilon >0$.
Suppose $m$ is a positive integer and $m'$ is a divisor of $m$ relatively prime to $D(F)$
that satisfies $(m')^{1+\epsilon}\ge m^{(2/d)+\epsilon}/|D(F)|^{1/d(d-1)}$. Then the number
of primitive solutions to (1) is at most
$$\left ( 5600d+{700\over \epsilon }\right )d^{\omega (m')}.$$
\endproclaim

The constants 5600 and 700 here carry no particular importance beyond specificity. The 
major improvement over the result of Bombieri and Schmidt is that the quantity
$\omega (m')$ is possibly much smaller than $\omega (m)$. 
In the same paper, Stewart explicitly constructed forms of various degrees to show
lower bounds for
the number of primitive solutions to (1). In so doing, he was lead to the following.

\proclaim{Conjecture (Stewart)} There is an absolute constant $c_0$ such that, for
all forms $F$ as in the theorem above, there is a positive bound $C$ (depending on $F$)
such that (1) has at most $c_0$ primitive solutions for all $m\ge C$.
\endproclaim

In this paper we will obtain results which simplify and strengthen Stewart's. Perhaps as
important is that our method
has the added benefit of providing good
heuristics for the conjecture above. In order to state our main results, 
we introduce a bit more notation.

Denote the set of places of $\bq$ by $M(\bq )$. For any $v\in M(\bq )$ we let
$|\cdot |_v$ denote the usual $v$-adic absolute value on $\bq$ and $\bq _v$ denote
the topological completion of $\bq$ with respect to this absolute value, though we will continue to
use $|\cdot |$ for the usual Euclidean absolute value. We fix
algebraic closures $\overline{\bq _v}$ for each of these and
assume that our original absolute values on $\bq$
are extended to the $\overline{\bq _v}$'s. As usual, we identify the finite places with positive
primes.

Any form $F(X,Y)\in\bq [X,Y]$ factors completely into a product of linear forms over
some splitting field:
$$F(X,Y)=\prod _{i=1}^d L_i(X,Y).$$
This splitting field may be embedded into any $\overline{\bq _v}$; we abuse notation
somewhat and write the above for the factorization of $F$ over $\overline{\bq _v}$ for
all places $v\in M(\bq )$. These linear factors are only unique up to a scalar multiple,
of course. We say a
linear factor $L_i(X,Y)$ is defined over $\bq _v$ if all possible quotients of coefficients 
are in $\bq _v$. 
For any form $F(X,Y)\in\bz [X,Y]$ 
and place $v\in M(\bq )$, set $c_F(v)$ to be the number of 
linear factors that are defined over $\bq _v$. For any integer $m>1$ set
$$c_F(m)=\prod \Sb p|m\\ p\ \text{prime}\endSb c_F(p).$$

\proclaim{Theorem 1} Let $F(X,Y)\in\bz [X,Y]$ be a homogeneous polynomial
of degree $d\ge 2$ with non-zero discriminant and content 1 and suppose
$m$ is a positive integer with $|m|_p<|D(F)|_p$ for all primes $p|m$. 
Then the primitive $(x,y)\in\bz ^2$ with $m|F(x,y)$ are contained in 
$c_F(m)$ sublattices of $\bz ^2$ of determinant $m.$
\endproclaim

In particular, there are no solutions to (1) if $c_F(m)=0$. In other words, $c_F(m)=0$ implies that there
is some local obstruction to solving (1).

Given a sublattice $\Lambda\subseteq\bz ^2$ of relatively large determinant, we can
provide an upper bound not just on the number of primitive solutions to (1),
but even to the related inequality
$$|F(x,y)|\le m.\tag 1'$$

\proclaim{Theorem 2} Let $F(X,Y)\in\bz [X,Y]$ be a homogeneous polynomial
of degree $d\ge 3$ with non-zero discriminant and content 1. Suppose
$m$ is a positive integer and $\Lambda\subseteq
\bz ^2$ is a sublattice with $\det (\Lambda )=Am^{2/d}/|D(F)|^{1/d(d-1)}$ for some $A>0$. If
$A\ge 5^4,$
then the number of primitive lattice points $(x,y)\in\Lambda$ that are solutions to
(1') is less than
$$2+2d\left ( 11+{31\over\log (d-1)}+{\log\left ( {2\log m\over d\log A}+2\right )\over\log (d-1)}\right ).$$
If $A<5^4,$
then the number of solutions is less than
$${2\cdot 5^4\over A}
\left (
2+2d\left ( 11+{31\over\log (d-1)}+{\log\left ( {\log m\over 2d\log 5}+2\right )\over\log (d-1)}\right )
\right ).$$
\endproclaim

\proclaim{Corollary} Let $F(X,Y)\in\bz [X,Y]$ be a homogeneous polynomial of degree $d\ge 3$ 
with non-zero 
discriminant and content 1. Suppose $m$ is a positive integer and $m'$ is a divisor of
$m$ relatively prime to $D(F)$ that satisfies $m'=Am^{2/d}/|D(F)|^{1/d(d-1)}$ for some $A\ge 1$. Then
the number of primitive solutions $(x,y)\in\bz ^2$ to (1') with $m'|F(x,y)$ is less than
$$2500d\left ( 43+{\log \big (2+ \log m/(1+\log A )\big )\over\log (d-1)}\right )c_F(m').$$ 
\endproclaim

\demo{Proof} One readily checks that for $d\ge 3,$
$(d/2)\log A\ge 1+\log A$ if $A\ge 5^4$ and $2d\log 5\ge 1+\log A$ if $1\le A<5^4$.
The rest follows from Theorems 1 and 2.
\enddemo

Comparing the Corollary with the result of Stewart above, our constants are in the same ballpark even though
we are estimating more than just the solutions to (1). Moreover, we can easily replace the 2500 with
2 once $A\ge 5^4$. Our real improvement is
replacing Stewart's $d^{\omega (m')}$ term with $c_F(m')$. Clearly $c_F(m')\le d^{\omega 
(m')}$ always, but will typically be much smaller. Also, the estimate in Stewart's result
tends to infinity as the divisor $m'$ approaches $m^{2/d}/|D(F)|^{1/d(d-1)}$, whereas our 
estimate is bounded above by a constant multiple of $c_F(m')d\log\log m /\log (d-1)$. 

Again, the main novelty of our approach is the lattices; primitive solutions to our Thue equation
(1) are elements of certain sublattices of $\bz ^2$. One may view Theorem 2 as attempting to
limit the number $c_F(m')$ of sublattices considered by allowing the divisor $m'$ to be small, yet still
allowing for a ``good" upper bound for the number of solutions (at least a result as strong as Schmidt's
conjecture above, say). Another approach is to take $m'$ as large as possible and see that there
are very few solutions to (1') in the associated lattices. This approach can be used to give good
heuristics for Stewart's conjecture above. We consider this approach now.

For a given form $F(X,Y)\in\bz [X,Y]$
with non-zero discriminant and positive integer $m$, we set $m(F)$ to be the largest divisor $m'$
of $m$ with $|m'|_p<|D(F)|_p$ for all primes $p|m'$. Thus, $m(F)$ is the largest divisor of $m$
satisfying the hypotheses of Theorem 1. Clearly any solution $(x,y)$ to (1) satisfies $m(F)|F(x,y)$.
For any form $F$ we write $\uF$ for the coefficient vector. Given an $F(X,Y)$ with its factorization
into linear forms as above and an $\epsilon >0$, we
say a non-zero $(x,y)\in\bz ^2$ is $\epsilon$-{\it exceptional} if
$${|L_i(x,y)L_j(x,y)|\over |\det (\uL _i^{tr},\uL _j^{tr})|}\le {1\over \Vert (x,y)\Vert ^{\epsilon}}$$
for some indices $i\neq j$, where $\uL ^{tr}$ denotes the transpose of $\uL$ and $\Vert\cdot\Vert$ 
denotes the usual Euclidean norm on $\bc ^2$. In other words, the $\epsilon$-exceptional points are
the points dealt with by Roth's theorem (or the Subspace theorem the way we have formulated things
here). As is well-known, the number of such exceptional points is bounded above by an explicit
function of $\epsilon$ and the degree $d$ of $F$, thus justifying the ``exceptional" moniker.

\proclaim{Theorem 3} Let $F(X,Y)\in\bz [X,Y]$ be a homogeneous polynomial
of degree $d\ge 5$ with non-zero discriminant and content 1. For any positive $\epsilon < d-4$
there is a positive $c(F,\epsilon )$ depending only
on $F$ and $\epsilon$ such that if $m\ge c(F,\epsilon )$ and $\Lambda$ is a sublattice of $\bz ^2$
with determinant $m(F)$, then there is at most one
pair of primitive solutions $\pm (x,y)\in\Lambda$ to (1') that is not $(d-4-\epsilon )$-exceptional. 
Further, any such  primitive solution must satisfy
$\Vert (x,y)\Vert < m(F)^{(1/2)-(\epsilon/3(2+\epsilon ))}$.
\endproclaim

A major point here is that for $m$ large enough, any primitive solution to (1) that isn't
exceptional in the sense of Thue-Siegel-Roth is a non-zero 
lattice point of length much
less than $\sqrt{\det (\Lambda )}$. But now lattices with such a point are exceptional in 
that their first successive minima (in the sense of Minkowski) is smaller than typical.	

\proclaim{Theorem 4} Let $\delta >0.$ 
The proportion of sublattices $\Lambda\subseteq\bz ^2$ of determinant $m$ 
with a primitive $(x,y)\in\Lambda$ satisfying $\Vert (x,y)\Vert \le m^{(1/2)-\delta }$ is
$O(m^{-2\delta })$ as $m\rightarrow\infty$, where the implicit constant is absolute.
\endproclaim

We believe that Theorems 3 and 4 in conjunction with Theorem 1 lend credence to Stewart's 
conjecture when the degree $d\ge 5$.
For these degrees, any primitive solutions to (1) either give exceptionally good 
rational approximations to the lines where the form $F$ vanishes or arise from lattices with
exceptionally small first minima. Therefore the existence of any primitive solutions to (1) is indeed
exceptional, at least for $m$ large enough. 

\head Lattices Arising from Thue Equations\endhead

Our goal in this section is to prove Theorem 1.
For any form $F(X,Y)\in\bz [X,Y]$ we write $\uF$ for the coefficient vector.
We first note that $(x,y)\in\bq ^2$
is a primitive integral point only if $\Vert (x,y)\Vert _p\le 1$ for all primes $p$, where $\Vert\cdot
\Vert _p$ denotes the supnorm. Also, given any
positive integer $m'$, $m'|F(x,y)$ if and only if $|F(x,y)|_p\le |m'|_p$ for all primes $p$.
We will use the following non-archimedean version of [9, Lemma 4].

\proclaim{Lemma 1} Let $K$ be a topologically complete field with respect to
a non-archimedean 
absolute value $|\cdot |$ and $L_1(\uX ),\ldots ,L_n(\uX )\in K[\uX ]$ be
$n$ linearly independent 
linear forms. Let $\| \cdot \|$ denote the supnorm on $K^n$. Suppose $\ux\in K^n$ 
and $j$ is such that 
$${|L_j(\ux )|\over \|\uL _j\|}\ge {|L_i(\ux )|\over \|\uL _i\|}$$
for $i=1,\ldots ,n$. Then
$${|L _j(\ux )|\over\|\uL _j\|}
\ge {\|\ux\| |\det (\uL ^{tr}_1,\ldots ,\uL ^{tr}_n)|
\over \prod _{i=1}^n\|\uL _i\|}.$$
\endproclaim

\demo{Proof} The statement is obvious is $\ux =\uo$, so suppose otherwise.
Then without loss of generality we may assume $\| \uL _i\| =1$
for all $i$ and $\|\ux\| =1$. Let $T$ denote the $n\times n$
matrix with rows $\uL _i$ and write
$$\gathered {\frak m}=\min \Sb \uy\in K^n\\ \| \uy \|  =1\endSb
\left \{\| T \uy ^{tr} \| \right \}\\
{\frak M}=\max \Sb \uy\in K^n\\ \|\uy\| =1\endSb
\left \{\|T \uy ^{tr}\| \right \}.\endgathered$$

Suppose $\| T\ux _1^{tr}\| ={\frak m}$ and $\|\ux _1\| =1.$ 
Choose $\ux _2,\ldots ,
\ux _n\in K^n$, all of length 1, that also satisfy
 $|\det (\ux _1^{tr}, \ldots
,\ux _n^{tr})| =1.$ We then have
$$\aligned |\det (T)|=|\det (T)||\det (\ux _1^{tr},\ldots ,\ux _n^{tr})| &=
|\det (T\ux _1^{tr},\ldots ,T\ux _n^{tr})| \\
&\le\prod _{l=1}^n\|T\ux _l^{tr}\| \\
&\le {\frak m}{\frak M}^{n-1}.\endaligned$$

Since $\| \uL _i\| =1$ for all $i$ and the absolute value is non-archimedean
we have ${\frak M}\le 1$, so that
${\frak m}\ge |\det (T)|.$ On the other hand, by our choice of $j$ we also have
$|L_j(\ux )| \ge |L_i(\ux )|$ for all $i=1,\ldots ,n$. Since $\|\cdot\|$ is the 
supnorm, these $n$ inequalities (and the definition of $T$) imply that
$|L_j(\ux )|\ge\|T\ux ^{tr}\|\ge {\frak m}.$ Thus
$$|L_j(\ux )|\ge {\frak m}\ge |\det (T)|=|\det (\uL _1^{tr},\ldots ,\uL _n^{tr})|.$$
\enddemo

\proclaim{Lemma 2} Let $F(X,Y)\in\bz [X,Y]$ be a homogeneous polynomial
with non-zero discriminant and content 1. 
Then for every primitive $(x,y)\in\bz ^2$ and every prime $p$ with $|F(x,y)|_p<|D(F)|_p$
there is a linear factor $L_p(X,Y)\in\bzp [X,Y]$ of $F$ with $\|\uL _p\| _p=1$
and $$|L_p(x,y)|_p\le |F(x,y)|_p,$$
with equality if $p\nmid D(F)$.
\endproclaim

\demo{Proof} Write
$$F(X,Y)=\prod _{i=1}^d M_i(X,Y),$$
where $M_i(X,Y)\in\bqpbar [X,Y]$ is a linear form for all $i=1,\ldots ,d$.
Let $|\cdot |_p$ be an absolute value on $\bqpbar$ that extends the usual
$p$-adic absolute value on $\bqp$.
Suppose $(x,y)\in\bz ^2$ is a primitive integral point and choose $i_0$ such that
$${|M_{i_0}(x,y)|_p\over\|\uM _{i_0}\| _p}=\min _{1\le i\le d}\left \{ 
{|M_i(x,y)|_p\over\|\uM _i\| _p}\right \}.$$
We note by Lemma 1 that
$${|M_i(x,y)|_p\over \|\uM _i\| _p}\ge {|\det (\uM _{i_0}^{tr},\uM _i^{tr})|_p\over
\|\uM _{i_0}\| _p\|\uM _i\| _p}\tag 2$$
for all $i\neq i_0$, since $\| (x,y)\| _p=1$. 

We claim that $M_{i_0}$ is defined over $\bqp$. Indeed, if this were not the 
case, then without loss of generality
there would be a $\sigma$ in the Galois group of $\bqpbar$ over $\bqp$ with
$\sigma (M_{i_0})=M_{i_1}$ for some $i_1\neq i_0$ between $1$ and $d$. We then
have $|M_{i_1}(x,y)|_p=|M_{i_0}(x,y)|_p$ and 
$\|\uM _{i_1}\|_p=\|\uM _{i_0}\|_p$ (see [1, chap. 2, Theorem 7]),
so that by Lemma 1 
$${|M_{i_0}(x,y)|_p\over \|\uM _{i_0}\| _p}\ge {|\det (\uM _{i_0}^{tr},
\uM _{i_1}^{tr})|_p\over\|\uM _{i_0}\| _p\|\uM _{i_1}\| _p}.\tag 3$$
Now by (2), (3), Hadamard's inequality and Gauss' lemma
$$\aligned
{|F(x,y)|_p\over \|\uF\| _p}={|F(x,y)|_p\over\|\uM _1\| _p\cdots\|\uM _d\| _p}&\ge
{|\det (\uM _{i_1}^{tr},\uM _{i_0}^{tr})|_p\over\|\uM _{i_0}\| _p\|\uM _{i_1}\| _p}
\prod _{i\neq i_0}{|\det (\uM _{i_0}^{tr},\uM _i^{tr})|_p\over
\|\uM _{i_0}\| _p\|\uM _i\| _p}\\
&\ge
\prod _{i\neq j}{|\det (\uM _i^{tr},\uM _j^{tr})|_p\over
\|\uM _i\| _p\|\uM _j\| _p}\\
&={|D(F)|_p\over \|\uM _1\| _p^{2(d-1)}\cdots \|\uM _d\| _p^{2(d-1)}}\\
&={|D(F)|_p\over \|\uF\| _p^{2(d-1)}}.\endaligned$$
Since the content of $F$ is 1, we get $|F(x,y)|_p\ge |D(F)|_p$ which contradicts our original hypothesis.
Thus $M_{i_0}$ is defined over $\bqp$.

Arguing exactly as above, but this time only using (2), we have
$${|F(x,y)|_p\over \|\uF\| _p}\ge {|M_{i_0}(x,y)|_p\over\|\uM _{i_0}\| _p}
{| D(F)|_p\over \|\uF\| _p^{2(d-1)}}.$$
This gives
$${|M_{i_0}(x,y)|_p\over\|\uM _{i_0}\| _p}\le {|F(x,y)|_p\over |D(F)|_p}.$$
Since $M_{i_0}$ is defined over $\bqp$, we let $L_p(X,Y)\in\bzp [X,Y]$ be
a scalar multiple of $M_{i_0}(X,Y)$ with $\|\uL _p\| _p=1$. We note that if
$p\nmid D(F)$, then Hadamard's inequality is actually an equality in all of
the above and Lemma 1 gives
$${|M_i(x,y)|_p\over \|\uM _i\| _p}=1$$
for all $i\neq i_0$. Thus
$$|L_p(x,y)|_p={|M_{i_0}(x,y)|_p\over\|\uM _{i_0}\| _p}=|F(x,y)|_p$$
if $p\nmid D(F)$.
\enddemo

\proclaim{Lemma 3} Suppose $p$ is a prime and $L(X,Y)\in\bzp [X,Y]$ is
a linear form with $\| \uL\| _p=1$. Let $\alpha _p$ denote the Haar measure
on $\bqp$ with $\alpha _p (\bzp )=1$. Then for all integers $c\ge 0$ the
set
$$S=\{ (x,y)\in\bzp ^2\: |L(x,y)|_p\le p^{-c}\}$$
is a $\bzp$-module with $\alpha _p^2(S)=p^{-c},$ where $\alpha _p^2$ denotes the
product measure on $\bq _p^2$.
\endproclaim

\demo{Proof} Clearly $S$ is a $\bzp$-module. Write
$L(X,Y)=aX+bY$ and set
$$M(X,Y)=\cases Y&\text{if $|a|_p=1$,}\\
X&\text{if $|a|_p<1$.}\endcases$$
Note that in the second case here we necessarily have $|b|_p=1$ since $\|\uL \| _p=1.$
In either case we easily have $\|\uM\| _p=1$ and $|\det (\uL ^{tr},\uM ^{tr})|_p=1$.
Now choose $\uz _1,\uz _2\in\bz ^2$ with
$$|L(\uz _1)|_p=|M(\uz _2)|_p=1,\qquad L(\uz _2)=M(\uz _1)=1.$$
Then
$$\bzp ^2=\{ x\uz _1+y\uz _2\: x,y\in \bz _p\},\qquad
S=\{ x\uz _1+y\uz _2\: x,y\in \bz _p,\ |x|_p\le p^{-c}\},$$
so that $\alpha ^2 _p(S)=p^{-c}\big ( \alpha _p(\bz _p)\big )^2\alpha _p^2(\bz _p^2)=p^{-c}.$
\enddemo

\proclaim{Lemma 4} Let ${\Cal S}$ be a finite set of prime numbers and
$L_p(X,Y)\in\bzp [X,Y]$ be a linear form with $\|\uL _p\| _p=1$ for all
$p\in {\Cal S}$. For each $p\in {\Cal S}$ let $a_p$ be a positive integer and set
$$S_p=\{ (u,v)\in\bzp ^2\: |L_p(u,v)|_p\le p^{-a_p}\}.$$
Set $S_p=\bzp ^2$ for all primes $p\not\in {\Cal S}$. 
Then 
$$\Lambda =\bigcap \Sb p\ \text{prime}\endSb \bq ^2\cap S_p$$
is a sublattice of $\bz ^2$ with
$$\det (\Lambda )=\prod _{p\in {\Cal S}} p^{a_p}.$$
\endproclaim

\demo{Proof} This follows immediately from Lemma 3 and general facts on lattices and
$\bz _p$-modules (see [11, chap. 3], for example).\enddemo

\demo{Proof of Theorem 1} 
If $(x,y)\in\bz ^2$ is a primitive point with $m|F(x,y)$ and $p$ is a prime dividing $m$, then
$|F(x,y)|_p\le |m|_p<|D(F)|_p$ so
by Lemma 2 there is a linear factor $L_p(x,y)\in\bz _p [X,Y]$ of $F$ with 
$\| \uL _p\| _p=1$
and $|L_p(x,y)|_p\le |F(x,y)|_p\le|m|_p.$
There are $c_F(p)$ possible linear factors here
by definition, whence $c_F(m)$ choices in total when we consider all primes dividing
$m$. Now suppose
for each prime $p|m$ we have chosen a linear factor $L_p(X,Y)\in\bz _p[X,Y]$ of $F$
with $\|\uL _p\| _p=1$. Then by Lemma 4 the set of all $(x,y)\in\bz ^2$ with
$|L_p(x,y)|_p\le |m|_p$ for all primes $p|m$ is a sublattice of $\bz ^2$ of
determinant $m$.
\enddemo

\head Proof of Theorem 2\endhead

If $F(X,Y)$ is a any form and
$\Lambda = \bz \uz_1\oplus\bz\uz _2$ is a lattice, then considering solutions $\uz\in
\Lambda$ to (1') is the same as considering solutions $(x,y)\in\bz ^2 $ to the
inequality $|\FL (x,y)|\le m$, where the form $\FL (X,Y):= 
F(X\uz _1+Y\uz _2)$. 
The choice of basis is not unique here of course.
We may also view $\FL (X,Y)$ as a composition $F\circ T$, where
$T\in\gltwo (\br )$ sends the canonical basis of $\bz ^2$ to a basis of $\Lambda$.
Note that a different choice of basis amounts to multiplying $T$ by an element of
$\gltwo (\bz )$. 

Our proof will involve various heights which we now define. For any form $F$ written
as a product of linear forms, $F(X,Y)= \prod _{i=1}^dL_i(X,Y)$, we set
$$\hofF =\prod _{i=1}^d\| \uL _i\|,$$
where $\|\cdot \|$ denotes the usual $L_2$ norm on $\bc ^2$. We set
$$\MofF =\min _{T\in\gltwo (\bz )}{\Cal H}(F\circ T)\qquad \text{and}\qquad
\mofF =\min \Sb T\in\gltwo (\br )\\ |\det (T)|=1\endSb {\Cal H}(F\circ T).$$
We remark that in general (see [10, Lemma 1]) for any form $F$ of degree $d$ 
and any $T\in\gltwo (\br ),$
$$\gathered
D(F\circ T)=D(F)\det (T)^{d(d-1)}\\
{\frak m}(F\circ T)=\mofF |\det (T)|^{d/2}\\
\MofF\ge\mofF\ge |D(F)|^{1/2(d-1)}.\endgathered\tag 4$$
In particular, we see that $|D(\FL )|$, $\mofFL$ and $\MofFL$ are all well-defined 
(i.e., are independent of the particular choice of basis) and satisfy
$$|D\big (\FL ) |= |D(F)|\det (\Lambda )^{d/2},\quad
\MofFL\ge\mofFL = \mofF \det (\Lambda )^{d/2}.\tag 4'$$
For a given positive integer $m$ we set
$\MofFLm$ to be the minimum of ${\Cal H}(\FL )$ over all bases $\uz _1,\ \uz _2$ of
$\Lambda$ with $\uz _1$ a solution to (1'), assuming such a primitive solution exists.

The main idea for determining solutions to (1') is
to say that some linear factor of $F$ must be relatively small for a given solution. 
For example, suppose we rewrite
$F(X,Y)=a\prod _{i=1}^d(X-\alpha _iY)$. Now if 
$(x,y)\in\bz ^2$ is any solution to (1'), then
$$|\alpha _i-x/y|\le {d2^{d-1}m\hofF ^{d-2}\over |y|^d|D(F)|^{1/2}}=d2^{d-1}(\hofF /m)^{d-2}
{m^{d-1}\over |D(F)|^{1/2}}{1\over |y|^d}\tag 5$$
for some index $i$ by [6 chap. 3, Lemmas 3A and 3B].
An alternative to (5) is that any solution $\ux$ to (1') satisfies
$${|L_i(\ux )L_j(\ux )|\over |\det (\uL _i^{tr},\uL _j^{tr})|}\le {2^{d-2}m\hofF ^{d-2}\over
\|\ux \| ^{d-2}|D(F)|^{1/2}}=(2\hofF /m)^{d-2}{m^{d-1}\over |D(F)|^{1/2}}{1\over \|\ux\| ^{d-2}}
\tag 5'$$
for some indices $i$ and $j$; this is [9, Lemma 5] (with the constants made
explicit). 

Considering either (5) or (5'), one can see that the major goal is to estimate
those solutions $\ux = (x,y)$ to (1') either with $|y|$ or $\| \ux \|$ ``small," so
that any remaining ``large" solutions are $\epsilon$-exceptional for some $\epsilon >0$.
Such ``large" solutions may be dealt with using gap arguments and ultimately a quantitative
version of Roth's theorem. 

We will use the following as the main part of our proof of Theorem 2.

\proclaim{Proposition} 
Suppose $F(X,Y)\in\bz [X,Y]$ is a form of degree $d\ge 3$ with
non-zero discriminant and content 1, $m$ is a positive integer and $\Lambda\subset
\br ^2$ is a lattice with $\det (\Lambda )=Am^{2/d}/|D(F)|^{1/d(d-1)}$. Then 
$\MofFLm \ge A^{d/2}m$ and if $A\ge 5^4$ the
number of primitive lattice points that are solutions to (1') is less than
$$2+2d\left ( 11+{\log 2^{10}3^35^3\over\log (d-1)}+{\log 2^93^35^2\over\log (d-5/4)}
+{\log\left ( {\log m\over \log (\MofFLm )-\log m}+2\right )\over\log (d-1)}\right ) .$$
\endproclaim

The proof of the Proposition will rely on a few lemmas, though we note that the
inequalities $\MofFLm\ge\MofFL \ge A^{d/2}m$ follow directly from the definitions,
(4), (4') and the hypotheses. To prove the Proposition we
obviously may assume there is a primitive lattice point $\uz _0\in\Lambda$ 
that is a solution to (1') since otherwise there is nothing to prove. Given this 
assumption, we choose a basis $\uz _0,\ \uz _0'$ of $\Lambda$ such that  $\uz _0$
is a solution to (1') and $\MofFLm = {\Cal H}(\FL )$. We will write  
$$F(X,Y)=\prod _{i=1}^d L_i(X,Y),\qquad \FL (X,Y)= \prod _{i=1}^dXL_i(\uz _0)+YL_i (
\uz _0')=F(\uz _0)\prod _{i=1}^d X+\alpha _iY,$$
where $\alpha _i=L_i(\uz _0')/L_i(\uz _0)$. 
For notational convenience, in what follows we will denote the quantity 
$\MofFLm /m$ by $B$. The hypothesis that $A\ge 5^4$ is thus equivalent to the assumption that $B\ge 5^{2d}$.

With the above conventions in place, we see by (4') and (5) that for any solution 
$\uz =x\uz _0+y\uz _0'\in\Lambda$ to (1') with $y\neq 0$ there is some index $i$ with
$$\aligned |\alpha _i -x/y|&\le
{d2^{d-1}m^{d-1}B^{d-2}\over |y|^d|D(\FL )|^{1/2}}\\
&= {d2^{d-1}m^{d-1}B^{d-2}\over |y|^d|D(F)|^{1/2}\det (\Lambda )^{d(d-1)/2}}\\
&\le {d2^{d-1}B^{d-2}\over |y|^d (5^4)^{d(d-1)/2}}\\
&<{B^{d-2}\over 2|y|^d}.\endaligned\tag 6$$
We may utilize a standard gap principle argument to estimate those solutions with 
$|y|> B$, for example (see Lemma 7 below). Eventually we come to the point where
a quantitative version of Roth's theorem is invoked  (Lemma 8). But before we do that,
we deal with those solutions where $|y|$ is smaller. The following is a variation on
[6, chap. 3, Lemma 5B]. 

\proclaim{Lemma 5} 
For every primitive lattice point $\uz =x\uz _0+y\uz _0'\in\Lambda$ with $y\neq 0$ that
is a solution to (1'),
there are $\psi _1(\uz ),\ldots ,\psi _d(\uz )\in [0,1]$ that, if not zero, are at least
$1/(2d)$, satisfy $\sum _{i=1}^d\psi _i (\uz )\ge 1/2$, and also 
$${|L_i(\uz _0)|\over |L_i(\uz )|}\ge \big ( B^{\psi _i (\uz )}-2\big ) |y|$$
for all $i=1,\ldots ,d.$
\endproclaim

\demo{Proof} 
We first claim that $2|L_{i_0}(\uz _0)|\le |L_{i_0}(\uz )|$ for some index $i_0$. 
Indeed, if this
were not the case then $\Lambda ':=\bz \uz _0\oplus\bz\uz $ is a sublattice of $\Lambda$
and $F_{\Lambda '}(X,Y):= F(X\uz _0+Y\uz )$ satisfies
$${\Cal H}(F_{\Lambda '}) ^2 =\prod _{i=1}^d|L_i(\uz _0)|^2+|L_i(\uz )|^2
<\prod _{i=1}^d5|L_i(\uz _0)|^2
\le 5^dm^2$$
since $\uz _0$ is a solution to (1'). But now by (4), (4') and the hypotheses
we have a contradiction:
$$\aligned 5^{d/2}m > {\Cal H}( F _{\Lambda '} )
&\ge \mofF \det (\Lambda ')^{d/2}\\
&\ge \mofF\det (\Lambda )^{d/2}\\
&\ge |D(F)|^{1/2(d-1)}\det (\Lambda )^{d/2}\\
&\ge 5^{2d}m.\endaligned$$

With the claim shown, choose an index $i_0$ with $2|L_{i_0}(\uz _0)|\le |L_{i_0}(\uz )|$.
Since $\uz$ is a primitive lattice point there is a $\uz '\in\Lambda$ with
$\Lambda =\bz\uz\oplus\bz\uz '$. Further, we may add any integer multiple of $\uz$ to
$\uz '$ here. Thus, we may choose $\uz '$ such that
$\alpha :=\Re\big ( L_{i_0}(\uz ')/L_{i_0}(\uz )\big )$ satisfies $|\alpha |\le 1/2.$
We now write $\uz _0=z\uz  +z'\uz '$ for some $z,z'\in\bz$
with $|z'|=[\Lambda\:\bz\uz _0\oplus\bz\uz ]=|y|$.
For any linear form $L(X,Y)$ we have
$${L (\uz _0)\over L(\uz )}=z+z'{L(\uz ')\over L(\uz )}.$$
In particular, using $L=L_{i_0}$ we see that $|z+z'\alpha |\le 1/2$, 
and for all $i=1,\ldots ,d$
$$\aligned {|L_i(\uz _0)|\over |L_i(\uz )|}&=\left | z+z'{L _i(\uz ')\over  L_i(\uz )}
\right |\\
&= \left | z'\left ({L _i(\uz ')\over  L_i(\uz )}-\alpha\right )+z+z'\alpha \right |\\
&\ge  \left | z'\left ({L _i(\uz ')\over  L_i(\uz )}-\alpha\right )\right |-|z+z'\alpha |\\
&\ge |z'|\left ( \left | {L _i(\uz ')\over  L_i(\uz )}\right |-{1\over 2}
\right )-{1\over 2}\\
&\ge |y| \left (\left | {L _i(\uz ')\over  L_i(\uz )}\right |+1-2\right ).
\endaligned\tag 7$$
Since $|F(\uz )|\le m$ and $\Lambda =\bz\uz\oplus\bz\uz '$,
$$\aligned \prod _{i=1}^d\left ( 1+\left | {L_i(\uz ')\over L_i(\uz )}\right |\right )&\ge
\prod _{i=1}^d\sqrt{ 1+\left | {L_i(\uz ')\over L_i(\uz )}\right |^2}\\
&={1\over |F(\uz )|}\prod _{i=1}^d\sqrt {|L_i(\uz )|^2+|L_i(\uz ')|^2}\\
&\ge {\MofFLm\over |F(\uz )|}\\
&\ge B.\endaligned\tag 8$$

We define $\psi _i(\uz )$ by
$$B^{\psi _i(\uz )}=\cases B&\text{if $\left | {L_i(\uz ')\over L_i(\uz )}\right |+1\ge B$,}\\
1&\text{if $\left | {L_i(\uz ')\over L_i(\uz )}\right |+1< B^{1/(2d)}$,}\\
\left | {L_i(\uz ')\over L_i(\uz )}\right |+1&\text{otherwise.}
\endcases$$
Now by construction $0\le \psi _i(\uz )\le 1$ for all $i=1,\ldots ,d$ and any
$\psi _j(\uz )\ge 1/{2d}$ if it isn't zero. 
We have $\sum _{i=1}^d\psi _i(\uz )\ge 1$
if any $\psi _j(\uz )=1$, so suppose $\psi _i(\uz )<1$ for
all $i=1,\ldots ,d$. Then by (8)
$$\aligned
B^{1/2}\prod _{i=1}^d B^{\psi _i(\uz )}&>
\prod \Sb 1\le i\le d\\ \psi _i(\uz )=0\endSb \left ( 1+\left | {L_i(\uz ')\over L_i(\uz )}
\right |\right )
\prod \Sb 1\le i\le d\\ \psi _i(\uz )>0\endSb B^{\psi _i(\uz )}\\
&=\prod \Sb 1\le i\le d\endSb \left ( 1+\left | {L_i(\uz ')\over L_i(\uz )}
\right |\right )\\
&\ge B.\endaligned$$
This shows that $\sum _{i=1}^d\psi _i(\uz )\ge 1/2$ in all cases. Also by construction
$B^{\psi _i(\uz )}\le \left | {L_i(\uz ')\over L_i(\uz )}\right | +1$ for all $i$, so
that the remaining desired inequalities follow from (7).
\enddemo

\proclaim{Lemma 6} 
For all $c>0$ there are less than $2d(2c+1)$ primitive solutions $\uz =x\uz _0+y\uz _0'
\in\Lambda$ to (1') with $y\neq 0$ and $|y|\le B^c.$
\endproclaim

We thus are able to rather efficiently estimate solutions where $|y|\le B^c$ for any
fixed constant $c$. In particular, though it's certainly possible to improve
upon particular aspects of Lemma 5, there wouldn't be much to gain (the exception being
if one could improve upon $B$, specifically, if one could replace $B$ by a larger
quantity in terms of $m$ or $F$). However, we remark that the hypothesis $\det (\Lambda )\ge 5^4m^{2/d}
|D(F)|^{1/d(d-1)}$ can be relaxed to $\det (\Lambda )\ge 5^4(m/\mofF )^{2/d}$, both in Lemma 5 and here in
Lemma 6.
 
\demo{Proof} By Lemma 5
$${|L_i(\uz )|\over |L_i(\uz _0)|}=|\alpha _i-x/y|\le {1\over (B^{\psi _i(\uz )}-2)|y|^2}
\tag 9$$
for all solutions $\uz =x\uz _0+y\uz _0'\in\Lambda$ to (1') with $y\neq 0$
and all $i=1,\ldots ,d$.

Let ${\Cal S}$ denote the set of primitive solutions $\uz =x\uz _0+x\uz _0'\in\Lambda$ to 
(1') with $1\le y\le B^c.$
For the moment fix an index $i$ and consider the sum $\sum \psi _i(\uz )$ over all 
$\uz\in {\Cal S}$. Obviously we may restrict to solutions with $\psi _i(\uz )\neq 0;$
we arrange these solutions $\uz _l=x_l\uz _0+y_l\uz _0',\ l=1,\ldots ,n$ so that
$y_l\le y_{l+1}$ for all $l$. Then by Lemma 5 and (9)
$$\aligned
{1\over |y_ly_{l+1}|}&\le \left | {x_l\over y_l}-{x_{l+1}\over y_{l+1}}\right |\\
&\le \left | \alpha _i-{x_l\over y_l}\right |
+ \left | \alpha _i-{x_{l+1}\over y_{l+1}}\right |\\
&\le  {1\over (B^{\psi _i(\uz _l)}-2)|y_l|^2}+
{1\over (B^{\psi _i(\uz _{l+1})}-2)|y_{l+1}|^2}\\
&\le  {1\over (B^{\psi _i(\uz _l)}-2)|y_l|^2}+
{1\over (B^{\psi _i(\uz _{l+1})}-2)|y_ly_{l+1}|},
\endaligned$$
whence
$$|y_{l+1}|\ge (B^{\psi _i(\uz _l)}-2)\big ( 1-(B^{\psi _i(\uz _{l+1})}-2)^{-1}\big )|y_l|.
\tag 10$$

Since $\psi _i(\uz _l)\ge {1/2d}$ for all our $\uz _l$ and $B\ge 5^{2d},$ we have
$B^{\psi _i(\uz _l)}\ge 5$ and thus
$B^{\psi _i(\uz _l)}-3\ge B^{\psi _i(\uz _l)\log 2/\log 5}$. We now 
repeatedly apply (10) to get 
$$\aligned B^c\ge |y_n|&\ge (B^{\psi _i(\uz _1)}-2)(B^{\psi _i(\uz _2)}-3)\cdots
(B^{\psi _i(\uz _{n-1})}-3)\big ( 1- (B^{\psi _i(\uz _n)}-2)^{-1})|y_1|\\
&>\prod _{l=1}^{n-1}(B^{\psi _i(\uz _l)}-3)\times \big ( 1- (1/3)\big )\\
&\ge (2/3)\prod _{l=1}^{n-1}B^{\psi _i(\uz _l)\log 2/\log 5}.\endaligned$$
Taking logarithms yields
$$c+\log _B (3/2)> \sum _{l=1}^{n-1}\psi _i(\uz _l)\log 2/\log 5,$$
and since $\psi _i(\uz _n)\le 1$,
$$c+\log _B (3/2)+\log 2/\log 5 > \sum _{l=1}^{n}\psi _i(\uz _l)\log 2/\log 5
=\sum \Sb \uz\in {\Cal S}\endSb \psi _i(\uz )\log 2/\log 5.$$
Finally, by Lemma 5 and this last inequality
$$\aligned 
|{\Cal S}|&\le\sum \Sb \uz\in {\Cal S}\endSb\sum _{i=1}^d 2\psi _i(\uz )\\
&=\sum _{i=1}^d\sum \Sb\uz\in {\Cal S}\endSb 2\psi _i(\uz )\\
&<\sum _{i=1}^d 2c+2\log _B (3/2)+2\log 2/\log 5\\
&\le \sum _{i=1}^d 2c+2\log (3/2)/\log (5^{2d})+2\log 2/\log 5\\
&= \sum _{i=1}^d 2c+\log (3/2)/d\log 5+2\log 2/\log 5\\
&\le \sum _{i=1}^d 2c+\log (3/2)/3\log 5+2\log 2/\log 5\\
&<d(2c+1).\endaligned$$

The same argument works for estimating the number of primitive solutions $x\uz _0+
y\uz _0'$ to (1') with
$1\le -y\le B^c$. 
\enddemo

\proclaim{Lemma 7} 
For all $C_2>C_1>B$, the number of primitive solutions $\uz =x\uz _0+y\uz _0'\in\Lambda$ 
to (1') with $C_1\le |y|\le C_2$ is less than
$$2d\left ( 1+{\log\big ( \log C_2/\log (C_1/B)\big )\over\log (d-1)}\right ).$$
\endproclaim

\demo{Proof} We will use (6). Suppose $x\uz _0+y\uz _0',\ x'\uz _0+y'\uz _0'\in\Lambda$
are primitive solutions to (1') with both
$$|\alpha _i -x/y|< {B^{d-2}\over 2|y|^d},\qquad
|\alpha _i -x'/y'|< {B^{d-2}\over 2|y'|^d}$$
for some index $i$. Suppose further that $y'\ge y>0.$ Then by the inequalities above
$$\aligned{1\over |yy'|}&\le \left | {x\over y}-{x'\over y'}\right |\\
&\le \left | \alpha _i-{x\over y}\right |+
\left | \alpha _i-{x'\over y'}\right |\\
&< {B^{d-2}\over 2|y|^d}
+{B^{d-2}\over 2|y'|^d}\\
&\le{B^{d-2}\over |y|^d},\endaligned$$
so that $|y'|\ge {|y|^{d-1}/B^{d-2}}.$
Hence if $x_1\uz _0+y_1\uz _0',\ x_2\uz _0+y_2\uz _0',\ldots $ are primitive solutions to 
(1') as above with
$C_1\le y_1\le y_2\le \cdots\le C_2,$ then
repeatedly applying the above inequality yields
$$C_2\ge y_{l+1}\ge {y_1^{(d-1)^l}\over B^{((d-1)^{l-1}+\cdots +1)(d-2)}}\ge
{C_1^{(d-1)^l}\over B^{(d-1)^l-1}}>(C_1/B)^{(d-1)^l}$$
for all $l\ge 1$. We take logarithms twice to get
$${\log\big ( \log C_2/(\log C_1/B)\big )\over\log (d-1)}>l.$$
Taking into account the $d$ possible indices $i$ and employing the same argument for
solutions with $y<0$ gives the lemma.
\enddemo

\proclaim{Lemma 8} 
Then there are less than
$$2d\left (4   +{\log 2^93^35^2\over\log (d-5/4)}\right )$$
primitive $\uz =x\uz _0+\uz _0'\in\Lambda$  solutions to (1') with 
$|y|\ge\max\{B^{4(d-1)},(8^d\MofFLm )^{2^{10}3^35^3}\}.$ 
\endproclaim

\demo{Proof} 
We note that the $\alpha _i$ are conjugate algebraic numbers with absolute height 
$$h(\alpha _i)^d={\Cal H}(\FL )=\MofFLm$$ 
(see [6, chap. 3, Lemma 2A], for example). Given a solution as in the lemma, by
(10) and the hypothesis $|y|\ge B^{4(d-1)}$ we have
$$\aligned
|\alpha _i-x/y|&\le {d2^{d-1}B^{d-2}\over |y|^d}\\
&<{B^{d-1}\over 2|y|^d}\\
&\le {1\over 2|y|^{d-1/4}}
\endaligned\tag 11$$
for some index $i$.
We claim that
$$|\alpha _i-x/y|<\big ( H(x/y)\big )^{-\sqrt{2d}(1+1/20)},\tag 12$$
where $H(x/y)=\sqrt{x^2+y^2}$ is the (absolute) height of $x/y$. 
To see this, we first note that $|x/y|<|\alpha _i| +1$, so that
$H(x/y)<(|\alpha _i|+2)|y|\le 3h(\alpha _i)^d|y|.$ Since $d\ge 3$, one readily
verifies that $d-1/4\ge\sqrt{2d}(1+1/10)$. Using this we easily get
$(3h(\alpha _i )^d)^{d-1/4}<(3h(\alpha _i)^d )^{\sqrt{d}}
<y^{\sqrt{2d}/20}$ (with quite a bit of room to spare, in fact). 
In addition, we also get
$$\aligned
y^{d-1/4}&\ge \left ( {H(x/y)\over 3h(\alpha _i)^d }\right )^{d-1/4}\\
&>H(x/y)^{\sqrt{2d}(1+1/10)}y^{-\sqrt{2d}/20}\\
&\ge H(x/y)^{\sqrt{2d}(1+1/10)}H(x/y)^{-\sqrt{2d}/20}\\
&= H(x/y)^{\sqrt{2d}(1+1/20)}.\endaligned$$
Therefore, (12) follows from (11). 

According to [6, chap. 2, Theorem 6] (with $m=2$ and $\chi = 1/20$ there), 
the rational solutions $x/y$ to (12) satisfy $H(x/y)\le 
(8h(\alpha _i))^{d2^{10}3^35^3}$ or $w\le H(x/y)<w^{2^93^35^2d^2}$ for some
$w>1$. The first option here is ruled out for us by hypothesis since $H(x/y)\ge |y|$. Hence
it remains to estimate the number of primitive solutions $(x,y)$ to (11) with
$w/(3h(\alpha _i)^d )\le |y|<w^{2^93^35^2d^2}$. We clearly may assume that 
$w\ge (8h(\alpha _i))^{4d}$. 

Suppose $(x_0,y_0), \ (x_1,y_1),\ldots$ are the
primitive solutions to (11) with $y_i>0$ and arranged so that $0<y_0\le y_1\le\cdots $
We then have
$$\aligned{1\over |y_ly_{l+1}|}&\le \left | {x_l\over y_l}-{x_{l+1}\over y_{l+1}}\right |\\
&\le \left | \alpha _i-{x_l\over y_l}\right |+
\left | \alpha _i-{x_{l+1}\over y_{l+1}}\right |\\
&< {1\over 2|y_l|^{d-1/4}}
+{1\over 2|y_{l+1}|^{d-1/4}}\\
&\le{1\over |y_l|^{d-1/4}},\endaligned$$
so that $|y_{l+1}|\ge |y_l|^{d-5/4}$ for all $l\ge 0$. Moreover,
since $w\ge (8^dh(\alpha _i)^d )^4$ and $d\ge 3$ we have
$$w^{d-2}\ge (8h(\alpha _i))^{4d}
>(3h(\alpha _i)^d )^{2(d-2)}B^{2(d-2)}\ge (3h(\alpha _i)^d B)^{d-1},$$
so that also by (11)
$$y_1\ge {y_0^{d-1}\over B^{d-1}}\ge \left ( {w\over 3h(\alpha _i)^d B}\right ) ^{d-1}
>{w^{d-1}\over w^{d-2}}=w.$$
We thus have $y_l\ge w^{(d-5/4)^l}$ for all $l\ge 1$.
Now since all $y_l<w^{2^93^35^2d^2}$, we must have
$$l<{\log 2^93^35^2d^2\over\log (d-5/4)}<3+{\log 2^93^35^2\over\log (d-5/4)}.$$ 

Considering the $d$ possible indices $i$ above and accounting for those solutions 
with $y<0$ in the same manner completes the proof.
\enddemo

\demo{Proof of the Proposition} We first set $c=2$ in Lemma 6 to see that the number
of primitive solutions $\uz =x\uz _0+y\uz _0'\in\Lambda$ to (1') with
$1\le |y|\le B^2$ is less than $10d$. Next we set $C_1=B^2$ and $C_2=B^{4(d-1)}$ in
Lemma 7 to see that the number of solutions with $B^2\le |y|\le B^{4(d-1)}$
is less than 
$$ 2d\left ( 1+{\log \big ( \log B^{4(d-1)}/\log B\big )\over\log (d-1)}\right )
=2d\big ( 2+\log 4/\log (d-1)\big ).$$
If on the other hand we set 
$C_2=\big ( 8^d\MofFLm\big )^{2^{10}3^35^3}$, then (recall $B\ge 5^{2d}>8^d$) the number of
solutions with $B^2\le |y|\le \big ( 8^d\MofFLm \big )^{2^{10}3^35^3}$ is less than
$$\multline 2d\left ( 1+{\log\big ( 2^{10}3^35^3\log (8^d\MofFLm )/\log B\big )\over
\log (d-1)}\right ) 
\\ 
<2d \left ( 1+ {\log \big ( 2^{10}3^35^3(1+\log (\MofFLm )/\log B )\big )\over\log (d-1)}
\right )\\
=
2d \left ( 1+ {\log \big ( 2^{10}3^35^3(2+\log m /\log B )\big )\over\log (d-1)}
\right )\\
=
2d \left ( 1+ {\log 2^{10}3^35^3\over\log (d-1)}+{\log (2+\log m /\log B )\over\log (d-1)}
\right ).\endmultline$$
Therefore the number of solutions with
$B^2\le |y|\le\max\{ B^{4(d-1)}, \big ( 8^d\MofFLm \big )^{2^{10}3^35^3}\}$ is less  than
$$2d \left ( 2+{\log 2^{10}3^35^3\over\log (d-1)}+{\log (2+\log m /\log B )\over\log (d-1)}
\right ).$$
Combining this with Lemma 8, the number of solutions with $y\neq 0$ is less than
$$10d+
2d \left ( 6+ {\log 2^{10}3^35^3\over\log (d-1)}+{\log 2^93^35^2\over\log (d-5/4)}
+{\log (2+\log m /\log B )\over\log (d-1)}\right ).$$
Of course, we also have the two solutions $\pm\uz _0$ as well, giving the Proposition.
\enddemo
 
\demo{Proof of Theorem 2} Suppose first that $A\ge 5^4$. We may assume that there is a 
primitive solution $(x,y)\in\Lambda$ to (1'). We apply the Proposition, noting that
$$\log (\MofFLm /m)\ge\log (A^{d/2}),\qquad
{\log (2^{10}3^35^3)\over\log (d-1)}+{\log (2^93^35^2)\over\log (d-5/4)}<{31\over\log (d-1)}$$
since $d\ge 3$.
For $A<5^4$ we use the Proposition in conjuction with Lemma 2C (and
Remark 2D) of [6, chap. 3] as follows. Let $p$ be any prime satisfying $(5^4/A)
\le p\le 2(5^4/A)-1$ and
let $F$ be a form as in the Proposition except that $A<5^4$. Then there are $p+1$ forms
$G$ with $|D(G)|=|D(F)|p^{d(d-1)}$ and any primitive integer
solution $(x,y)$ to (1') is a primitive integral solution to $|G(x,y)|\le m$ for one of these
forms $G$. Since 
$$\det (\Lambda )={Am^{2/d}\over |D(F)|^{1/d(d-1)}}={Am^{2/d}p\over |D(G)|^{1/d(d-1)}}\ge {5^4m^{2/d}
\over |D(G)|^{1/d(d-1)}},$$
we may apply the Proposition to these $p+1\le 2(5^4/A)$ forms $G$ to prove the case of
Theorem 2 when $A<5^4$. 
\enddemo

\head Proof of Theorems 3 and 4\endhead

\demo{Proof of Theorem 3} Let $\Lambda\subseteq\bz ^2$ be a sublattice with $\det (\Lambda )=m(F)$.
Denote the successive minima of $\Lambda$ (with respect to the unit disk) by $\lambda _1\le\lambda _2$.
By Minkowski's theorem, 
$$\lambda _2^2\ge\lambda _1\lambda _2\ge (2^2/2!){\det (\Lambda )\over \pi}={2m(F)\over\pi}.\tag 13$$
We clearly have $m(F)\ge m/|D(F)|$. Thus for $m$ sufficiently large
(depending on both $F$ and $\epsilon$), 
$$m(F)^{1+\epsilon /12}> m,\qquad
m(F)^{\epsilon /12}\ge {2^{d-2}\hofF ^{d-2}\over (2/\pi)^{1+\epsilon /2}|D(F)|^{1/2}}.\tag 14$$ 
Now suppose $\ux =(x,y)\in\Lambda$ is a primitive lattice point with $\Vert (x,y)\Vert \ge \lambda _2$.
If $(x,y)$ is a solution to (1'), then by (5'), (13) and (14)
there  are indices $i\neq j$ with 
$$\aligned
{|L_i(\ux )L_j(\ux )|\over |\det (\uL _i^{tr},\uL _j^{tr}|}&\le {2^{d-2}m\hofF ^{d-2}\over\Vert \ux\Vert
^{d-2}|D(F)|^{1/2}}\\
&\le {2^{d-2}m\hofF ^{d-2}\over\lambda _2^{2+\epsilon}\Vert\ux\Vert ^{d-4-\epsilon} |D(F)|^{1/2}}\\
&\le {2^{d-2}m\hofF ^{d-2}\over m(F)^{1+\epsilon /2}(2/\pi )^{1+\epsilon /2}\Vert\ux\Vert
^{d-4-\epsilon}|D(F)|^{1/2}}\\
&={m\over m(F)^{1+\epsilon /4}}{2^{d-2}\hofF ^{d-2}\over m(F)^{\epsilon /4}(2/\pi )^{1+\epsilon /2}|D(F)|
^{1/2}}{1\over\Vert\ux\Vert ^{d-4-\epsilon }}\\
&<{1\over\Vert\ux \Vert ^{d-4-\epsilon }}.\endaligned$$
Thus $\ux =(x,y)$ is $(d-4-\epsilon )$-exceptional. For notational convenience temporarily set
$\delta = \epsilon /3(2+\epsilon )$. Now if $\ux =(x,y)\in\Lambda$
is a primitive solution to (1') with $\Vert\ux\Vert =\lambda _1$ and $\lambda _1\ge m(F)^{(1/2)-\delta},$
then as above (this time using the full strength of (14)) there are indices
$i\neq j$ with
$$\aligned
{|L_i(\ux )L_j(\ux )|\over |\det (\uL _i^{tr},\uL _j^{tr}|}&\le {2^{d-2}m\hofF ^{d-2}\over\Vert \ux\Vert
^{d-2}|D(F)|^{1/2}}\\
&= {2^{d-2}m\hofF ^{d-2}\over\lambda _1^{2+\epsilon}\Vert\ux\Vert ^{d-4-\epsilon} |D(F)|^{1/2}}\\
&< {2^{d-2}m\hofF ^{d-2}\over m(F)^{((1/2)-\delta )(2+\epsilon )}\Vert\ux\Vert
^{d-4-\epsilon}|D(F)|^{1/2}}\\
&= {2^{d-2}m\hofF ^{d-2}\over m(F)^{1+\epsilon /6}\Vert\ux\Vert
^{d-4-\epsilon}|D(F)|^{1/2}}\\
&={m\over m(F)^{1+\epsilon /12}}{2^{d-2}\hofF ^{d-2}\over m(F)^{\epsilon /12}|D(F)|^{1/2}}
{1\over\Vert\ux\Vert ^{d-4-\epsilon}}\\
&<
{1\over\Vert\ux\Vert ^{d-4-\epsilon}}.\endaligned$$
Thus $\ux$ is again $(d-4-\epsilon)$-exceptional. This shows that any primitive solution
$\ux =(x,y)\in\Lambda$ to (1') that is not $(d-4-\epsilon )$-exceptional must satisfy
$\lambda _2>\Vert\ux\Vert =\lambda _1<m(F)^{(1/2)-\delta}.$ There can be at most one pair $\pm (x,y)$
of such primitive lattice points by the definition of successive minima.
\enddemo

\demo{Proof of Theorem 4} The number of sublattices $\Lambda\subseteq\bz ^2$
with determinant $m$ is equal to $\sum _{n|m}n$ (see [5, \S 3], for example), thus the
number $N$ of such lattices satisfies $m\le N\ll m\log\log m$ (though we will only use the lower bound
here). On the other hand, any primitive $(x,y)\in\bz ^2$ is a lattice point in exactly one sublattice
$\Lambda\subseteq\bz ^2$ with $m=\det (\Lambda )>(\pi /2)\Vert (x,y)\Vert ^2$. Indeed, for such
a lattice we must have $\Vert (x,y)\Vert =\lambda _1<\lambda _2$ by (13), so that $\Lambda =\bz (x,y)
\oplus\bz (x',y')$ for some $(x',y')\in\bz ^2$ with $xy'-x'y=m.$
Since $(x,y)$ is a primitive point, there is an $(x',y')\in\bz ^2$ with $xy'-x'y=m.$ Moreover, from
elementary number theory any other such point is of the form $(x',y')+n(x,y)$ for some 
integer $n$. Hence any sublattice of determinant $m$ containing $(x,y)$ has the same basis, 
so there is only
one such sublattice. 

Now suppose $\delta >0.$ If $m^{2\delta}<\pi /2$ then there is nothing to prove since 
$O(m^{-2\delta})=O(1)$. Otherwise we have $(\pi /2)m^{1-2\delta}\le m$. 
By what we have shown, the number of sublattices $\Lambda\subseteq\bz ^2$ with $\det (\Lambda )=m$ that 
contain a primitive $(x,y)$ with $\Vert (x,y)\Vert\le m^{(1/2)-\delta }$
is equal to the number $N'$ of such primitive points. Clearly $N'$ is no greater than the total
number of integral points in the disk with radius $m^{(1/2)-\delta}$, which in turn is no greater than
$4\pi m^{1-2\delta}$. (The number of integral points in the disk of radius $r\ge 1$ is no more than
$\pi (r+1)^2\le 4\pi r^2$.) Thus
$N'\le  4\pi m^{1-2\delta}$ and so the proportion of such lattices $N'/N=O(m^{-2\delta }).$
\enddemo

\enddocument